\documentclass[12pt]{amsart}

 \newtheorem{thm}{Theorem}[section]
 \newtheorem{cor}[thm]{Corollary}
 \newtheorem{lem}[thm]{Lemma}
 \newtheorem{prop}[thm]{Proposition}
 \theoremstyle{definition}
 \newtheorem{defn}[thm]{Definition}
 \theoremstyle{remark}
 
 \numberwithin{equation}{section}

 \newcommand{\norm}[1]{\left\Vert#1\right\Vert}
 
 \newcommand{\C}{\mathbb{C}}
 
\begin{document}

\title[$C_0$-Hilbert Modules]
 {$C_0$-Hilbert Modules}

\author{Yun-Su Kim.  }


\email{kimys@indiana.edu}


\dedicatory{}



\newpage

\begin{abstract}We provide the definition and fundamental properties of algebraic
elements with respect to an operator satisfying hypothesis $(h)$.
Furthermore, we analyze Hilbert modules using $C_0$-operators
relative to a bounded finitely connected region $\Omega$ in the
complex plane.
\end{abstract}

\maketitle

\section*{Introduction}
The theory of contractions of class $C_{0}$ was developed by
Sz.-Nagy-Foias \cite{S2}, Moore-Nordgren \cite{M}, and
Bercovici-Voiculescu [2,3], and J.A. Ball introduced the class of
$C_{0}$-operators relative to a bounded finitely connected region
$\Omega$ in the complex plane, whose boundary $\partial\Omega$
consists of a finite number of disjoint, analytic, simple closed
curves. The theory of Hilbert modules over function algebras has
been developed by Ronald G. Douglas and Vern I. Paulsen \cite{D}.

We analyze Hilbert modules using $C_0$-operators relative to
$\Omega$. Every operator $T$ defined on a Hilbert space $H$
satisfying hypothesis $(h)$ is not a $C_0$-operator relative to
$\Omega$. Thus, we provide the definition of an \emph{algebraic
element with respect to} $T$.

If $B$ is the set of algebraic elements with respect to $T$, and
it is closed, then naturally we have a bounded operator $T_B$ from
the quotient space $H/B$ to $H/B$. In section 2, we discuss the
relationships between the algebraic elements with respect to $T_B$
in $H/B$ and the algebraic elements with respect to $T$ in $H$.

In section 3, we define a module action on a Hilbert space $H$ by
using a $C_0$-operator $T$ relative to $\Omega$, and introduce a
$C_0$-Hilbert module $H_T$. Naturally, this raises the following
question : \vskip0.1cm
 If every element of $H_T$ is algebraic with respect to
$T$ over $A$, then $T$ is either a $C_0$-operator or not.
\vskip0.1cm

In this paper, we consider a case in which the rank of the
$C_0$-Hilbert module $H_T$ is finite, and we show that if a
generating set $\{h_{1},\cdot\cdot\cdot,h_{k}\}(k<\infty)$ of a
Hilbert module $H_T$ over $A$ is contained in ${B}$, then $T$ is a
$C_0$-operator.

Furthermore, if $B$ is closed, then by using the Jordan model of a
$C_0$-operator $T$ relative to $\Omega$, we show that there are
locally maximal $C_{0}$-submodules
$M_{i}(i=0,1,2,\cdot\cdot\cdot)$ of $H_T$ such that
$M_{0}\subset{M_{1}}\subset{M_{2}}\subset\cdot\cdot\cdot$.
\vskip0.1cm

 The author would like to express her appreciation to
Professors Hari Bercovici and Ronald G. Douglas for making some
helpful comments on this paper.



\section{Preliminaries and Notation}\label{13}

\subsection{Hilbert Modules}
Let $X$ be a compact, separable, metric space and let $C(X)$
denote the algebra of all continuous complex-valued functions on
$X$. A \emph{function algebra} on $X$ is a closed subalgebra of
$C(X)$, which contains the constant functions and separates points
of $X$.

\begin{defn}
Let $F$ be a function algebra, and let $H$ be a Hilbert space. We
say that $H$ is a \emph{Hilbert module over} $F$ if there is a
separately continuous mapping $\phi:F\times{H}\rightarrow{H}$ in
each variable satisfying : \vskip0.2cm

(a) $\phi(1,h)=h,$

(b) $\phi(fg,h)=\phi(f,\phi(g,h)),$

(c) $\phi(f+g,h)=\phi(f,h)+\phi(g,h),$

(d)
$\phi(f,\alpha{h}+\beta{k})=\alpha{\phi(f,h)}+\beta{\phi(f,k)},$

\noindent for every $f$, $g$ in $F$, $h$, $k$ in $H$, and $\alpha$
$\beta$ in $\C$.
\end{defn}

We will denote $\phi(f,h)$ by $f.h$. For $f$ in $F$, we let
$T_{f}:H\rightarrow{H}$ denote the linear map $T_{f}(h)=f.h$. If
$H$ is a Hilbert module over $F$, then by the continuity in the
second variable we have that $T_{f}$ is bounded.

\begin{defn}
Let $H$ be a Hilbert module over $F$. Then the \emph{module bound
of }$H$, is
\begin{center}$K_{F}(H)=\inf\{K:\norm{T_{f}}\leq{K\norm{f}}\texttt{ for all }f\texttt{ in }A\}.$
\end{center}
We call $H$\emph{ contractive }if $K_{F}(H)\leq{1}$.
\end{defn}


If $H$ is a Hilbert module over $A$, then a set
$\{h_{\delta}\}_{\delta\in{\Gamma}}\subset{H}$ is called a
\emph{generating set} for $H$ if finite linear sums of the form
\begin{center}$\sum_{i}f_{i}.h_{\delta_{i}}, f_{i}\in{A},
\delta_{i}\in\Gamma$
\end{center}
are dense in $H$.
\begin{defn}
If $H$ is a Hilbert module over $A$, then rank$_{A}(H)$, the
\emph{rank of} $H$ \emph{over} $A$, is the minimum cardinality of
a generating set for $H$.
\end{defn}

In the last few decades, the theory of Hilbert modules over
function algebras has been developed by Ronald G. Douglas and Vern
I. Paulsen \cite{D}.
\subsection{A Functional Calculus.}\label{11}

Let $H$ be a Hilbert space.
 Recall that $H^{\infty}$
is the Banach space of all (complex-valued) bounded analytic
functions on the open unit disk $\textbf{D}$ with supremum norm
\cite{S2}. A contraction $T$ in $L(H)$ is said to be
\emph{completely nonunitary }if there is no invariant subspace $K$
for $T$ such that $T|K$ is a unitary operator.

B. Sz.-Nagy and C. Foias introduced an important functional
calculus for completely non-unitary contractions.
\begin{prop}\label{12}Let $T\in{L(H)}$ be a completely non-unitary
contraction. Then there is a unique algebra representation
$\Phi_{T}$ from $H^{\infty}$ into $L(H)$ such that :\vskip0.2cm

(i) $\Phi_{T}(1)=I_{H}$, where $I_{H}\in{L(H)}$ is the identity
operator;

(ii) $\Phi_{T}(g)=T$, if $g(z)=z$ for all $z\in\textbf{D}$;

(iii) $\Phi_{T}$ is continuous when $H^{\infty}$ and $L(H)$ are
given the weak$^\ast$-

\quad\quad topology.

(iv) $\Phi_{T}$ is contractive, i.e.
$\norm{\Phi_{T}(u)}\leq\norm{u}$ for all
$u\in{H^{\infty}}$.\end{prop}

We simply denote by $u(T)$ the operator $\Phi_{T}(u)$.

B. Sz.- Nagy and C. Foias \cite{S2} defined the \emph{class
$C_{0}$} relative to the open unit disk \textbf{D} consisting of
completely non-unitary contractions $T$ on $H$ such that the
kernel of $\Phi_{T}$ is not trivial.  If $T\in{L(H)}$ is an
operator of class $C_{0}$, then \begin{center}ker
$\Phi_{T}=\{u\in{H^{\infty}}:u(T)=0\}$\end{center} is a
weak$^{*}$-closed ideal of $H^{\infty}$, and hence there is an
inner function generating ker $\Phi_{T}$. The \emph{minimal
function} $m_{T}$ of an operator of class $C_{0}$ is the generator
of ker $\Phi_{T}$. Also, $m_{T}$ is uniquely determined up to a
constant scalar factor of absolute value one \cite{B1}. The theory
of class $C_{0}$ relative to the open unit disk has been developed
by B.Sz.- Nagy, C. Foias (\cite{S2}) and H. Bercovici (\cite{B1}).

\subsection{Hardy spaces}
We refer to \cite{16} for basic facts about Hardy space, and
recall here the basic definitions.
\begin{defn}
The space ${\: H^{2}(\Omega)}$ is defined to be the space of
analytic functions $f$ on $\Omega$ such that the subharmonic
function $|f|^{2}$ has a harmonic majorant on $\Omega$. For a
fixed $z_{0}$ $\in\Omega$, there is a norm on
$H^{2}$$($$\Omega$$)$ defined by

 \begin{center}     $\|f\|=\inf\{ u(z_{0})^{1/2}$: $u$ is a harmonic majorant of
      $|f|^{2}\}$.\end{center}
\end{defn}

 Let $m$ be harmonic measure for the point $z_{0}$, let
 $L^{2}(\partial{\Omega})$ be the $L^{2}$-space of complex valued
 functions on the boundary of $\Omega$ defined with respect to $m$,
 and let $H^{2}(\partial{\Omega})$ be the set of
 functions $f$ in $L^{2}(\partial{\Omega})$ such that
 $\int_{\partial{\Omega}} f(z)g(z) dz$ = 0 for every $g$ that is
 analytic in a neighborhood of the closure of $\Omega$.
 If $f$ is in $H^{2}(\Omega)$, then there is a function $f^{\ast}$ in
 $H^{2}(\partial{\Omega})$ such that $f({z})$ approaches
 $f^{\ast}(\lambda_{0})$ as $z$ approaches $\lambda_{0}$
 nontangentially, for almost every $\lambda_{0}$ relative to $m$. The map
 $f\rightarrow{f^{\ast}}$ is an isometry from $H^{2}(\Omega)$ onto
 $H^{2}(\partial{\Omega})$. In this way, $H^{2}(\Omega)$ can be
 viewed as a closed subspace of $L^{2}(\partial{\Omega})$.

 A function $f$ defined on $\Omega$ is in $H^{\infty}(\Omega)$ if
 it is holomorphic and bounded. $H^{\infty}(\Omega)$ is a closed
 subspace of $L^{\infty}({\Omega})$ and it is a Banach
 algebra if endowed with the supremum norm. Finally, the mapping
 $f\rightarrow{f^{\ast}}$ is an isometry of $H^{\infty}(\Omega)$
 onto a week$^{*}$-closed subalgebra of
 $L^{\infty}(\partial{\Omega})$.

\subsection{$C_{0}$-operators relative to $\Omega$}

We will present in this section the definition of
$C_{0}$-operators relative to $\Omega$. Reference to this material
is found in Zucchi \cite{20}.

Let ${H}$ be a Hilbert space and ${K_{1}}$ be a compact subset of
the complex plane. If $T$$\in$$L(H)$ and
$\sigma(T)$$\subseteq$$K_{1}$, for $r=p/q$ a rational function
with poles off $K_{1}$, we can define an operator $r(T)$ by
$q(T)^{-1}p(T)$.
\begin{defn}
If $T$$\in$$L(H)$ and $\sigma(T)$$\subseteq$$K_{1}$, we say that
$K_{1}$ is a \emph{spectral set} for the operator $T$ if
$\norm{r(T)}$$\leq$$\max\{$$|$$r(z)$$|$$:$ $z$$\in$$K_{1}\}$,
whenever $r$ is a rational function with poles off $K_{1}$.
\end{defn}
If $T$ $\in$$L(H)$ is an operator with $\overline{\Omega}$ as a
spectral set and with no normal summand with spectrum in
$\partial\Omega$, i.e., $T$ has no reducing subspace
$M$$\subseteq$$H$ such that $T|M$ is normal and
$\sigma(T|M)$$\subseteq$$\partial\Omega$, then we say that $T$
satisfies \emph{hypothesis (h)}.
\begin{prop} ([20], Theorem 3.1.4)
Let $T\in{L(H)}$ be an operator satisfying ${\:hypothesis}$ $(h)$.
Then there is a unique algebra representation $\Psi_{T}$ from
$H^{\infty}(\Omega)$ into $L(H)$ such that :

$(i) \Psi_{T}(1)$=$I_{H}$, where $I_{H}$$\in$$L(H)$ is the
identity operator;

$(ii) \Psi_{T}(g)$=$T$, where $g(z)$=$z$ for all $z$$\in$$\Omega$;

$(iii) \Psi_{T}$ is continuous when $H^{\infty}(\Omega)$ and
$L(H)$ are given the $weak^{*}$-

\quad\quad topology.

$(iv)$ $\Psi_{T}$ is contractive, i.e.,
$\norm{\Psi_{T}(f)}$$\leq$$\norm{f}$ for all
$f$$\in$$H^{\infty}(\Omega)$.
\end{prop}

From now on we will indicate $\Psi_{T}(f)$ by $f(T)$ for all
$f$$\in$$H^{\infty}(\Omega)$.

\begin{defn}\label{3}
An operator $T$ satisfying hypothesis (h) is said to be of
\emph{class} $C_{0}$ \emph{relative to} $\Omega$ if there exists
$u$ $\in$ $H^{\infty}(\Omega)$$\setminus$$\{{0}\}$ such that
$u(T)$=$0$.
\end{defn}

\section{Algebraic Elements with respect to an Operator satisfying hypothesis $(h)$}
Every operator $T$ satisfying hypothesis $(h)$ is not a
$C_{0}$-operator relative to $\Omega$, and so we provide the
following definition.

\begin{defn}
Let $T\in{L(H)}$ be an operator satisfying hypothesis (h). An
element $h$ of $H$ is said to be \emph{algebraic with respect to}
$T$ provided that $\theta(T)h=0$ for some
$\theta\in{H^{\infty}(\Omega)}\setminus\{0\}$.

If not, $h$ is said to be \emph{transcendental with respect to}
$T$.
\end{defn}
If $A$ is a closed subspace of $H$ generated by
$\{a_{i}\in{H}:i=1,2,3,\cdot\cdot\cdot\}$, then $A$ will be
denoted by $\bigvee_{n=1}^{\infty}a_{i}.$

\begin{prop}\label{6}
Let $T\in{L(H)}$ be an operator satisfying hypothesis (h).
\vskip.1cm

(a) If $h\in{H}$ is algebraic with respect to $T$, then so is any
element

\quad in $\bigvee_{n=0}^{\infty}T^{n}h.$ \vskip.1cm

(b) If $h\in{H}$ is transcendental with respect to $T$, then so is
$T^{n}h$ for

\quad any $n=0,1,2,\cdot\cdot\cdot.$
\end{prop}
\begin{proof}
(a) Let $\theta\in{H^\infty(\Omega)}\setminus\{0\}$ such that
$\theta(T)h=0$.

Then for any $n=0,1,2,\cdot\cdot\cdot,$
\[\theta(T)(T^{n}h)=T^{n}(\theta(T)h)=0.\]

It follows that $\theta(T)h^{\prime}=0$ for any
$h^{\prime}\in\bigvee_{n=0}^{\infty}T^{n}h.$ \vskip0.2cm

\noindent(b) Suppose that $T^{k}h$ is algebraic with respect to
$T$ for some $k>0$. Thus there is
$f\in{H^{\infty}(\Omega)}\setminus\{0\}$ such that $f(T)T^{k}h=0.$

Let $f_{1}(z)=z^{k}f(z)$ for $z\in{\textbf{D}}$. Then
$f_{1}\in{H^{\infty}}(\Omega)\setminus\{0\}$ and
\[f_{1}(T)h=T^{k}{f(T)}h=f(T)T^{k}h=0\]
which contradicts to the fact that $h$ is transcendental with
respect to $T$.
\end{proof}
Note that $T^{0}$ denote the identity operator on $H$. \vskip0.1cm
By Theorem 1 in \cite{15}, if $h\in{H}$ is algebraic with respect
to $T$, then there is an inner function
$m_{h}\in{H^{\infty}(\Omega)}$ such that $m_{h}(T)h=0$ and $m_{h}$
is said to be a \emph{minimal function} of $h$ with respect to
$T$.

\begin{thm}\label{5}
Let $T\in{L(H)}$ be an operator satisfying hypothesis (h), and
$B=\{h\in{H}:h$ is algebraic with respect to $T\}$. \vskip0.2cm

(a) If $M=\{h_{i}:i=1,2,\cdot\cdot\cdot,k\}(k<\infty)$ is
contained in $B$, then so

\quad is $\bigvee_{n=0}^{\infty}T^{n}M$. \vskip0.2cm

(b) $B$ is a subspace of $H$.
\end{thm}
\begin{proof}
(a) By Proposition \ref{6} (a),
\begin{equation}\label{22}m_{h_{i}}(T)(T^{n}h_{i})=0\end{equation} for any
$i=1,\cdot\cdot\cdot,k$ and $n=0,1,2,\cdot\cdot\cdot$.

Let $\theta=m_{h_{1}}\cdot\cdot\cdot{m_{h_{k}}}$. Then
$\theta\in{H^{\infty}}(\Omega)\setminus\{0\}$, and
$\theta=\theta_{i}m_{h_{i}}$ for some
$\theta_{i}\in{H^{\infty}(\Omega)\setminus\{0\}}$. Thus, by
equation (\ref{22}),
\begin{equation}\label{21}\theta(T)(T^{n}h_{i})=\theta_{i}(T)m_{h_{i}}(T)(T^{n}h_{i})=0
\end{equation}
for any $i=1,\cdot\cdot\cdot,k$ and $n=0,1,2,\cdot\cdot\cdot$.

If $x\in{\bigvee_{n=0}^{\infty}T^{n}M}$, then there is a sequence
$\{x_{n}\}_{n=1}^\infty$ such that
\begin{center}$\lim_{n\rightarrow\infty}x_{n}=x$ and
$x_{n}=\sum_{i=1}^{k}a_{n,i}P_{n,i}(T)h_{i}$\end{center} for some
$a_{n,i}\in\C$ and a polynomial $P_{n,i}$. Then, equation
(\ref{21}) implies that
\begin{center}$\theta(T)(x_{n})=\theta(T)(\sum_{i=1}^{k}a_{n,i}P_{n,i}(T)h_{i})=\sum_{i=1}^{k}a_{n,i}\theta(T)P_{n,i}(T)h_{i}=0.$
\end{center}
It follows that $\theta(T)(x)=0$ for any
$x\in\bigvee_{n=0}^{\infty}T^{n}M$. Thus $x\in{B}$.\vskip0.2cm

(b) Clearly, $0\in{B}$. For $h_{1}$ and $h_{2}$ in $B$, if
$m_{1}(T)h_{1}=m_{2}(T)h_{2}=0$, where
$m_{i}(i=1,2)\in{H^{\infty}}(\Omega)\setminus\{0\}$, then
\[(m_{1}m_{2})(T)(\alpha_{1}h_{1}+\alpha_{2}h_{2})=0\]
for any $\alpha_{i}(i=1,2)$ in $\C$. Thus $B$ is a subspace of
$H$.

\end{proof}

Note that $B$ does not need to be closed.

If $T$ is a bounded operator on $H$ and $M$ is a (closed)
invariant subspace for $T$, then we can define a bounded operator
${T_{M}}:H/M\rightarrow{H/M}$ defined by
\[{T_{M}}([h])=[Th]\]
where $H/M$ is the quotient space. Since $M$ is $T$-invariant,
${T_{M}}$ is well-defined. Clearly, ${T_{M}}$ is a bounded
operator on $H/M$.

Let $R(\Omega)$ be the algebra of rational functions with poles
off $\overline{\Omega}$. We will say that a (closed) subspace $N$
is $R(\Omega)$-\emph{invariant} (or \emph{rationally invariant})
for an operator $T$ if it is invariant under $u(T)$ for every
$u\in{R(\Omega)}$.

If $N$ is a $R(\Omega)$-invariant subspace for an operator $T$
satisfying hypothesis $(h)$, then we can define
$\theta({T_{N}}):H/N\rightarrow{H/N}$ by
\[\theta({T_{N}})([h])=[\theta(T)h]\]
for $\theta\in{H^{\infty}(\Omega)}$ and $[h]\in{H/N}$. Since $N$
is $R(\Omega)$-invariant for the operator $T$, ${T_{N}}$ is
well-defined. Clearly, ${T_{N}}$ is a bounded operator on $H/N$.

\begin{defn}
Let $T\in{L(H)}$ be an operator satisfying hypothesis $(h)$ and
$M$ be an invariant subspace for $T$. An element $[h]$ of $H/M$ is
said to be \emph{algebraic with respect to} ${T_{M}}$ provided
that $\theta({T_{M}})[h]=0$ for some
$\theta\in{H^{\infty}(\Omega)}\setminus\{0\}$.

If not, $h$ is said to be \emph{transcendental with respect to}
${T_{M}}$.
\end{defn}
\begin{prop}
Let $T\in{L(H)}$ be an operator satisfying hypothesis (h) and
$B=\{h\in{H}:h$ is algebraic with respect to $T\}$. If $B$ is
closed, then it is $R(\Omega)$-invariant.

\end{prop}
\begin{proof}
Let $h\in{B}$ and $u\in{R(\Omega)}$. Then there is a nonzero
function $\phi\in{H^{\infty}(\Omega)}$ such that $\phi(T)h=0$.

It follows that $u(T)\phi(T)h=\phi(T)(u(T)h)=0$, that is,
$u(T)h\in{B}$.
\end{proof}

\begin{thm}\label{9}
Let $T\in{L(H)}$ be an operator satisfying hypothesis (h) and
$B=\{h\in{H}:h$ is algebraic with respect to $T\}$. If $B$ is a
closed subspace of $H$, then the following statements are
equivalent: \vskip0.2cm

(i) $[a]\in{H/B}$ is algebraic with respect to ${T_{B}}$.

(ii) $a$ is algebraic with respect to $T$.

\end{thm}
\begin{proof}
$(i)\rightarrow{(ii)}$ Since $[a]\in{H/B}$ is algebraic with
respect to ${T_{B}}$, there is a nonzero function $\theta_{1}$ in
$H^{\infty}(\Omega)$ such that $\theta_{1}(T)a\in{B}$.

It follows that
\begin{equation}\label{10}\theta_{2}(T)(\theta_{1}(T)a)=0\end{equation}
for some $\theta_{2}\in{H^{\infty}(\Omega)}\setminus\{0\}$.

Let
$\theta_{3}=\theta_{1}\cdot\theta_{2}\in{H^{\infty}(\Omega)}\setminus\{0\}$.
Then by equation (\ref{10}), $\theta_{3}(T)a=0$, and so $a\in{B}$.
\vskip0.2cm $(ii)\rightarrow{(i)}$ If $a\in{H}$ is algebraic with
respect to ${T}$, then there is a nonzero function $\theta$ in
$H^{\infty}(\Omega)$ such that $\theta(T)a=0$. Since $0\in{B}$,
$\theta(T_{B})[a]=[\theta(T)a]=0$.
\end{proof}

\begin{cor}\label{14}
Under the same assumption as Theorem \ref{9}, the following
statements are equivalent: \vskip0.2cm

(i) $[a]\in{H/B}$ is algebraic with respect to ${T_{B}}$.

(ii) $[a]=[0]$.
\end{cor}
\begin{proof}
By Theorem \ref{9}, it is clear.

\end{proof}

\begin{cor}

Let $T\in{L(H)}$ be an operator satisfying hypothesis (h) and
$M\subset{B}$ is a $R(\Omega)$-invariant subspace for $T$. Then
the following statements are equivalent: \vskip0.2cm

(i) $[a]\in{H/M}$ is algebraic with respect to ${T_{M}}$.

(ii) $a$ is algebraic with respect to $T$.

\end{cor}
\begin{proof}
It can be proven by the same way as the proof of Theorem \ref{9}.
\end{proof}

\begin{cor}

Let $T\in{L(H)}$ be an operator satisfying hypothesis (h) and
$M\subset{B}$ is a $R(\Omega)$-invariant subspace for $T$. Then
the following statements are equivalent: \vskip0.2cm

(i) $[a]\in{H/M}$ is transcendental with respect to ${T_{M}}$.

(ii) $a$ is transcendental with respect to $T$.

\end{cor}

We recall that if $K$ is a Hilbert space, $H$ is a subspace of
$K$, $V\in{L(K)}$, and $T\in{L(H)}$, then $V$ is said to be a
\emph{dilation} of $T$ provided that
\begin{equation}\label{11}T=P_{H}V|H.\end{equation} If $T$ and $V$
are operators satisfying hypothesis $(h)$ and $V$ is a
$C_{0}$-operator relative to $\Omega$ satisfying equation
(\ref{11}), then $V$ is said to be a $C_{0}$\emph{-dilation} of
$T$. We will not discuss about $C_0$-dilation any more in this
paper.

\begin{lem}
Let $T\in{L(H)}$ be an operator satisfying hypothesis (h) and
$B^{\prime}=\{h\in{H}:h$ is transcendental with respect to $T\}$.
It $h\in{B^{\prime}}$, then $u(T)h\in{B^\prime}$ for any
$u\in{R(\Omega)\setminus\{0\}}$
\end{lem}
\begin{proof}
Suppose that there is an element $h$ in $B^\prime$ such that
$u(T)h$ is algebraic with respect to $T$ for some
$u\in{R(\Omega)}\setminus\{0\}$.

Thus there is a nonzero function $\phi\in{H^{\infty}(\Omega)}$
such that \begin{equation}\label{15}\phi(T)u(T)h=0.\end{equation}

Let $\theta=\phi\cdot{u}$. Then
$\theta\in{H^{\infty}(\Omega)\setminus\{0\}}$ such that
$\theta(T)h=0$ by equation (\ref{15}). It contradicts to the fact
that $h\in{B^{\prime}}$.

\end{proof}

\section{$C_{0}$-Hilbert Modules}

Let $H$ be a Hilbert space and $F$ be a function algebra on $X$.
Then $H$ is a Hilbert module over $F$ with the module action
$F\times{H}\rightarrow{H}$ given by
\[f.h=f(x)h\]
for a fixed $x\in{X}$. Let $H_x$ denote this Hilbert module over
$F$. Clearly, $H_x$ is a contractive Hilbert module over $F$ for
any $x\in{X}$. \vskip0.2cm Similarly, for an operator $T$ on $H$
satisfying hypothesis $(h)$, if $A\subset{H^{\infty}}(\Omega)$ is
a function algebra over $\overline{\Omega}$ such that every
polynomial is contained in $A$, then $H$ is a Hilbert module over
$A$ with the module action $A\times{H}\rightarrow{H}$ given by
\begin{equation}\label{7}f.h=f(T)h.\end{equation} In this paper, $H_T$ denotes this Hilbert module
over $A\subset{H^{\infty}}(\Omega)$. Clearly, $H_T$ is a
contractive Hilbert module over $A$.

In this section, $A$ denotes a function algebra over
$\overline{\Omega}$ such that every polynomial is contained in $A$
and $A\subset{H^{\infty}}(\Omega)$.
\begin{defn}
If $T\in{L(H)}$ is a $C_{0}$-operator relative to $\Omega$, then
$H_{T}$ is called a $C_0$\emph{-Hilbert module}.

\end{defn}

\begin{defn}
Let $H$ and $K$ be Hilbert modules over $A$. Then a \emph{module
map} $X:H\rightarrow{K}$ is a bounded, linear map satisfying
$X(f.h)=f.(Xh)$ for all $f$ in $A$, and $h$ in $H$. Two Hilbert
modules are \emph{similar} if there is an invertible module map
from $H$ onto $K$, and are said to be \emph{isomorphic} if there
is a module map from $H$ onto $K$ which is a unitary.

\end{defn}

\begin{prop}\label{1}
For operators $T_{i}(i=1,2)$ in $L(H)$ satisfying hypothesis (h),
if $T_1$ and $T_2$ are similar operators, then $H_{T_1}$ and
$H_{T_2}$ are similar Hilbert modules over $A$.
\end{prop}
\begin{proof}
Let a module map $G:H\rightarrow{H}$ denote the similarity such
that $GT_{1}=T_{2}G$.

Define a linear map $Y:H_{T_{1}}\rightarrow{H_{T_{2}}}$ by
\begin{equation}\label{16}Y(f.h)=f.(Gh)\end{equation} for $f\in{A}$ and $h\in{H_{T_{1}}}$.

Let $f_{1}.h_{1}=f_{2}.h_{2}$ for $f_{i}\in{A}$ and
$h_{i}\in{H_{T_{1}}}$. Then
\begin{equation}\label{8}f_{1}(T_{1})h_{1}-f_{2}(T_{1})h_{2}=0.\end{equation} Since
$GT_{1}=T_{2}G$, equation (\ref{8}) implies that
\[f_{1}(T_{2})Gh_{1}=Gf_{1}(T_{1})h_{1}=Gf_{2}(T_{1})h_{2}=
f_{2}(T_{2})Gh_{2}.\] It follows that
$f_{1}.(Gh_{1})=f_{2}.(Gh_{2})$, that is, $Y$ is well-defined.

For $h\in{H_{T_{1}}}$,
\begin{equation}\label{17}Y(h)=Y(1.h)=G(h).\end{equation} By equations (\ref{16}) and (\ref{17}),
we can conclude that $Y$ is a module map. Since $G$ is bijective,
so is $Y$.

\end{proof}

\begin{cor}
For operators $T_{i}(i=1,2)$ in $L(H)$ satisfying hypothesis (h),
if $T_1$ and $T_2$ are unitarily equivalent, then $H_{T_1}$ and
$H_{T_2}$ are isomorphic.

\end{cor}

\begin{proof}

It is proven by the same way as the proof of Proposition \ref{1}.

\end{proof}
If $T\in{L(H)}$ is an operator satisfying hypothesis $(h)$ and $M$
is a submodule of $H_{T}$ over $A$, then by the definition of
module action given in equation (\ref{7}), we have that $M$ is
$T$-invariant. Furthermore, $M$ is a invariant subspace for each
operator $u(T)$ where $u\in{A}$.

\begin{defn}
Let $T\in{L(H)}$ be an operator satisfying hypothesis $(h)$. If
$M$ is a submodule of $H_{T}$ (over $A$) such that
$T|{M}:M\rightarrow{M}$ is a $C_{0}$-operator relative to
$\Omega$, then $M$ is said to be a $C_{0}$-\emph{submodule} (over
$A$) of $H_{T}$

\end{defn}
\begin{defn}
Let $T\in{L(H)}$. If there is an element $h\in{H}$ which is not in
the kernel of $T$ such that $\{T^{n}h:n=0,1,2,\cdot\cdot\cdot\}$
is not linearly independent, then $T$ is said to be
\emph{dependent}.

\end{defn}
\begin{thm}\label{25}
If $T\in{L(H)}$ is a dependent operator satisfying hypothesis (h),
then $H_{T}$ always has a nonzero $C_{0}$-submodule $M$.
\end{thm}
\begin{proof}
Since $T$ is dependent, there is a nonzero element $h$ in $H$ such
that $\{T^{n}h:n=0,1,2,\cdot\cdot\cdot\}(\neq\{0\})$ is linearly
dependent. It follows that
\[\sum_{n=0}^{k}a_{n}T^{i_{n}}h=0,\]
for some nonzero polynomial
$p(z)=\sum_{n=0}^{k}a_{n}z^{i_{n}}(z\in{\textbf{D}}).$

Let $M$ be the closed subspace of $H$ generated by
$\{\theta(T)h:\theta\in{A}\}$ and
$M^{\prime}=\{f\in{A}:f(T)h=0\}$. Since $p\in{M^{\prime}}$,
$M^{\prime}$ is not empty.

Clearly, $f.k$ is in $M$ for every $f$ in $A$ and $k$ in $M$ and
so $M$ is a submodule of $H_T$.

For any $\theta\in{A}$ and $f\in{M^\prime}$,
\[f(T)\theta(T)h=\theta(T)f(T)h=0.\]
It follows that $f(T)h^{\prime}=0$ for any $f\in{M^\prime}$ and
$h^{\prime}\in{M}$.

\noindent Therefore, $T_{0}=T|{M}$ is a $C_{0}$-operator relative
to $\Omega$, and so $M$
 is a $C_{0}$-submodule of $H_T$.
\end{proof}

\begin{defn}
Let $T\in{L(H)}$ be an operator satisfying hypothesis $(h)$. A
$C_{0}$-submodule $M$ of $H_T$ over $A$ is said to be
\emph{maximal} provided that there is no submodule $M^{\prime}$ of
$H_T$ over $A$ such that $M\subset{M^{\prime}}$ and $T|M^\prime$
is a $C_{0}$-operator relative to $\Omega$.

\end{defn}

\begin{cor}
Let $T\in{L(H)}$ be an operator satisfying hypothesis (h).

If $M$ is a maximal $C_{0}$-submodule of $H_{T}$ and
$h\in{H_{T}\setminus{M}}$, then
$\{T^{n}h:n=0,1,2,\cdot\cdot\cdot\}$ is linearly independent.
\end{cor}
\begin{proof}
Suppose that there is an element $h\in{H_{T}\setminus{M}}$ such
that $\{T^{n}h:n=0,1,2,\cdot\cdot\cdot\}$ is linearly dependent.

If $M^{\prime}$ is the closed subspace of $H_{T}$ generated by
$\{\theta(T)h:\theta\in{A}\}$, then by Theorem \ref{25},
$T|{M^\prime}$ is a $C_{0}$-operator relative to $\Omega$. Since
$T|M$ and $T|{M^\prime}$ are $C_{0}$-operators relative to
$\Omega$, there are nonzero functions
$\theta_{i}\in{H^{\infty}(\Omega)}(i=1,2,)$ such that
\begin{equation}\label{26}
\theta_{1}(T|M)=0\texttt{ and }\theta_{2}(T|{M^\prime})=0.
\end{equation}
It follows that $\theta_{1}\theta_{2}(T|M\vee{M^\prime})=0$, that
is, $T|M\vee{M^\prime}$ is also a $C_0$-operator relative to
$\Omega$. By maximality of $M$, $M\vee{M^\prime}=M$ which
contradicts to the fact that $h\in{M^{\prime}\setminus{M}}$.

\end{proof}

For an operator satisfying hypothesis $(h)$, $T\in{L(H)}$,
$h\in{H}$ is said to be \emph{algebraic with respect to $T$ over}
$A$, provided that
\begin{center}$\theta(T)h=0$ for some $\theta\in{A\setminus\{0\}}$.\end{center}

If $B=\{h\in{H}:h$ is algebraic with respect to $T$ over $A\}$,
then we could raise the question of whether the following sentence
is true or not; \vskip0.1cm

If every element of $H_T$ is algebraic with respect to $T$ over
$A$, then $T$ is a $C_0$-operator. \vskip0.1cm

In the next Theorem, we provide a condition in which that sentence
is true.

\begin{thm}\label{20}
Let $T\in{L(H)}$ be an operator satisfying hypothesis $(h)$. If
$H_T$ is a Hilbert module over $A$ with a generating set
$\{h_{1},\cdot\cdot\cdot,h_{k}\}(k<\infty)$ and $h_{i}\in{B}$ for
$i=1,2,\cdot\cdot\cdot,k$, then $H_{T}=B$ and $T$ is a
$C_0$-operator.

\end{thm}
\begin{proof}
Since $h_{i}\in{B}$, there is a nonzero function $m_i$ in $A$ such
that $m_{i}(T)h_{i}=0$ for $i=1,2,\cdot\cdot\cdot,k$. Then for any
$f\in{A}$,
$m_{i}(T)(f.h_{i})=m_{i}(T)f(T)h_{i}=f(T)m_{i}(T)h_{i}=0$. It
follows that $f.h_{i}\in{B}$ for any $f\in{A}$.

By Theorem \ref{5} $(b)$,
$\{\sum_{i=1}^{k}f_{i}.h_{i}:f_{i}\in{A}\}$ is contained in $B$.
Since $\{\sum_{i=1}^{k}f_{i}.h_{i}:f_{i}\in{A}\}$ is dense in
$H_T$, it is enough to prove that $B$ is a closed subspace of
$H_T$.

Let $b$ be an element in the closure of $B$ in the norm topology
induced by the inner product defined in $H_T$. Then, there is a
sequence $\{b_{n}\}_{n=1}^\infty$ in
$\{\sum_{i=1}^{k}f_{i}.h_{i}:f_{i}\in{A}\}$ such that
$\lim_{n\rightarrow\infty}b_{n}=b$.

Define a function $m=m_{1}\cdot\cdot\cdot{m_k}$. Then, for any
$f_{i}\in{A}$,
\begin{equation}\label{18}
m(T)(\sum_{i=1}^{k}f_{i}.h_{i})=m(T)(\sum_{i=1}^{k}f_{i}(T)h_{i})=\sum_{i=1}^{k}f_{i}(T)m(T)h_{i}=0.
\end{equation}
Equation (\ref{18}) implies that $m(T)(b_{n})=0$ for any
$n=1,2,\cdot\cdot\cdot$. Thus $m(T)b=0$ so that $b\in{B}$.
Therefore, $H_{T}=B$.

Since $m(T)b=0$ for any $b\in{B}(=H_T)$, $m(T)=0$ which proves
that $T$ is a $C_0$-operator.

\end{proof}
Recall that a nonzero function $\theta$ in $H(\Omega)$ is said to
be \emph{inner} if $|\theta|$ is constant almost everywhere on
each component of $\partial\Omega$. Then the \emph{Jordan block}
$S(\theta)$ is an operator acting on the space
$H(\theta)=H^{2}(\Omega)\ominus\theta{H^{2}(\Omega)}$ as follows :
\begin{center}$S(\theta)=P_{H(\theta)}S|{H(\theta)},$\end{center}
where $S\in{L(H^{2}(\Omega))}$ is defined by $(Sf)(z)=zf(z)$.

An operator $T\in{L(H)}$ is called a \emph{quasiaffine transform}
of an operator $T^{\prime}\in{L(H^{\prime})}(T\prec{T}^{\prime})$
if there exists an injective operator $X\in{L(H,H^{\prime})}$ with
dense range such that $T^{\prime}X=XT$. $T$ and $T^\prime$ are
\emph{quasisimilar} if $T\prec{T}^{\prime}$ and
$T^{\prime}\prec{T}$. \vskip0.2cm
\begin{prop}\cite{20}\label{2}
 Let $H$ be a separable Hilbert space and $T\in{L(H)}$ be an
operator of class $C_{0}$ relative to $\Omega$. Then there is a
family
$\{\theta_{i}\in{H^{\infty}(\Omega)}:i=0,1,2,\cdot\cdot\cdot\}$ of
inner functions such that \vskip0.2cm

(i) For $i=1,2,\cdot\cdot\cdot$, $\theta_{i}$ divides
$\theta_{i-1}$, that is, $\theta_{i-1}=\theta_{i}\varphi$ for some

\quad $\varphi\in{H^{\infty}(\Omega)}$.\vskip0.2cm

(ii) $T$ is quasisimilar to
$\bigoplus_{i=0}^{\infty}S(\theta_{i})$.
\end{prop}
If $T\in{L(H)}$ is a $C_{0}$-operator relative to $\Omega$, then
by Definition \ref{3}, ker $\Psi_{T}\neq\{0\}$ and there is an
inner function $\theta$, called a \emph{minimal function} of $T$,
in $H^{\infty}(\Omega)$ such that ker
$\Psi_{T}=\theta{H^{\infty}(\Omega)}$ \cite{20}. We denote by
$m_{T}$ the minimal function of $T$.

\begin{defn}
Let $M$ be a $C_{0}$-submodule of $H_{T}$ with the following
property ;\vskip0.2cm

If $M_{1}$ is a $C_{0}$-submodule of $H_T$ such that
$M\subset{M_{1}}$ and $m_{T|M_{1}}=m_{T|M}$, then
$M=M_1$.\vskip0.2cm

Then $M$ is said to be a \emph{locally maximal}
$C_{0}$\emph{-submodule} of $H_{T}$.

\end{defn}
\begin{thm}\label{19}
Let $H$ be a separable Hilbert space and $T\in{L(H)}$ be an
operator satisfying hypothesis $(h)$. If $B=\{h\in{H}:h$ is
algebraic with respect to $T$ over $A\}$ is a closed subspace of
$H$ and rank$_{A}H_{T}<\infty$, then there are locally maximal
$C_{0}$-submodules $M_{i}(i=0,1,2,\cdot\cdot\cdot)$ of $H_T$ such
that
\begin{center}$M_{0}\subset{M_{1}}\subset{M_{2}}\subset\cdot\cdot\cdot$.\end{center}
\end{thm}
\begin{proof}
Let $T^{\prime}=T|B$. For given element $h\in{B}$, we have a
function $m_{h}\in{A\setminus\{0\}}$ such that
\begin{center}$m_{h}(T)h=0.$
\end{center}
Then
$m_{h}(T)(\varphi.h)=m_{h}(T)\varphi(T)h=\varphi(T)m_{h}(T)h=0$
for any $\varphi$ in $A$. Thus, $B$ is a submodule of $H_{T}$ so
that $B=H_{T^\prime}$.

Since \begin{center}rank$_{A}H_{T^\prime}$=rank$_{A}B\leq$
rank$_{A}H_{T}<\infty$,\end{center} and every elements $h$ in $B$
is algebraic with respect to $T^\prime$ over $A$, Theorem \ref{20}
implies that $T^{\prime}=T|B$ is a $C_{0}$-operator.

Thus by Proposition \ref{2}, there are inner functions
$\theta_{i}(i=0,1,2,\cdot\cdot\cdot)$ such that $\theta_{i+1}$
divides $\theta_{i}$ and $T|B$ is quasisimilar to
$\bigoplus_{i=0}^{\infty}S(\theta_{i})$.

For each $\theta_{i}(i=0,1,2,\cdot\cdot\cdot)$, we have a bounded
linear operator $\theta_{i}(T):H\rightarrow{H}$ such that
\begin{center} $\theta_{i}(T)(f.h)=\theta_{i}(T)f(T)h=f(T)\theta_{i}(T)h=f.(\theta_{i}(T)h)$
\end{center}
for any $f\in{A}$ and $h\in{H}$. Thus
$\theta_{i}(T)(i=0,1,2,\cdot\cdot\cdot)$ is a module map.

It follows that $M_{i}=\ker(\theta_{i}(T))$ is a submodule of
$H_{T}$ and clearly, $T_{i}=T|M_{i}$ is a $C_{0}$-operator such
that $\theta_{i}(T_{i})=0$. Thus $M_{i}$ is a $C_{0}$-submodule of
$H_{T}$.

Let $i\in\{0,1,2,\cdot\cdot\cdot\}$ be given and $M$ be a
$C_{0}$-submodule of $H_T$ such that
\begin{equation}\label{23}M_{i}\subset{M}\texttt{ and }m_{T|M}=m_{T|M_i}.
\end{equation}
Since $m_{T|M_i}=\theta_i$, by equation (\ref{23}),
$m_{T|M}=\theta_i$. Thus, $\theta_{i}(T|M)=0$ so that
\begin{equation}\label{24}M\subset{\texttt{ ker}(\theta_{i}(T))=M_{i}}.
\end{equation}
From equations (\ref{23}) and (\ref{24}), $M=M_i$. Thus, $M_i$ is
a locally maximal $C_{0}$-submodule of $H_T$ for each
$i=1,2,3,\cdot\cdot\cdot$.

Since $\theta_{i+1}$ divides $\theta_{i}$ for
$i=0,1,2,\cdot\cdot\cdot$, $M_{i}\subset{M_{i+1}}$.

\end{proof}

In fact, in the proof of Theorem \ref{19}, $T|B$ is quasisimilar
to $\bigoplus_{i=0}^{k}S(\theta_{i})$ where $k\leq$
rank$_{A}H_{T}<\infty$. Thus, we have a finite number of locally
maximal $C_{0}$-submodules $M_{i}(i=0,1,2,\cdot\cdot\cdot,k)$.
\vskip0.2cm Naturally, the following question remains :  When is
$B$ closed?

\noindent However, we will not discuss this question in this
paper.

------------------------------------------------------------------------

\bibliographystyle{amsplain}
\bibliography{xbib}
\end{document}